\documentclass[dvipdfmx]{article}

\usepackage{amsmath,amsfonts,amssymb,amsthm,pb-diagram,picinpar,graphicx,color}
\usepackage{algorithm}
\usepackage{algorithmic}
\usepackage{colortbl}
\usepackage{enumerate}
\usepackage{hyperref}
\usepackage{bm}
\usepackage{tikz}
\usepackage{pgf}
\usepackage[T1]{fontenc}
\usepackage{enumitem}
\usetikzlibrary{shapes,matrix,calc,arrows,decorations.markings,knots,patterns,intersections,cd}

\newcommand{\CC}{\mathbb{C}}
\newcommand{\ZZ}{\mathbb{Z}}
\newcommand{\PP}{\mathbb{P}}
\newcommand{\mcC}{\mathcal{C}}
\newcommand{\mcB}{\mathcal{B}}
\newcommand{\bbL}{\mathbb {L}}

\newcommand{\fd}{\mathfrak{d}}

\newcommand{\pic}{\mathop{\mathrm{Pic}}\nolimits}
\newcommand{\amari}{\mathop{\mathrm{Rem}}\nolimits}
\newcommand{\Div}{\mathop{\mathrm{Div}}\nolimits}
\newcommand{\divi}{\mathop{\mathrm{div}}\nolimits}

\newcommand{\sr}{\mathop{\mathrm{sr}}\nolimits}

\newcommand{\NF}{\mathop{\mathrm{NF}}\nolimits}

\newtheorem{thm}{Theorem}[section]

\newtheorem{prop}[thm]{Proposition}
\newtheorem{prop0}{Proposition}

\theoremstyle{definition}
\newtheorem{defin}[thm]{Definition}
\newtheorem{defin0}[prop0]{Definition}

\theoremstyle{remark}
\newtheorem{rem}[thm]{Remark}

\renewcommand{\thesubparagraph}{\theparagraph.\@arabic\c@subparagraph}
\begin{document}
%
%
%
%


\begin{center}
{\Large \bf 
An explicit construction for  $n$-contact curves to a smooth cubic  via divisons and Zariski tuples
}

\bigskip

 {\large \bf Ai Takahashi and Hiro-o Tokunaga\footnote{Partially supported by Grant-in-Aid for Scientific Research C (20K03561)}}
 
 \end{center}
 
\begin{abstract}
Let $E$ be a smooth cubic. A plane curve $D$ is said to be an $n$-contact curve to $E$ if the intersection multiplicities at each intersection point between $E$ and $D$ is $n$.
In this note, we give an algorithm to produce $n$-contact curves to $E$ and consider its application.
\end{abstract}

\section*{Introduction}

In this note, all varieties are
 defined  over the field of complex numbers, $\CC$.
Let $C$ be a plane curve in $\PP^2$. We say that a  plane curve $D$ without multiple components is a contact curve to $C$ if  (i) $C\cap D$ consists of 
smooth points of $C$ and $D$ and (ii) for each $p \in C\cap D$, the intersection multiplicity  at $p$ is $\ge 2$.
In the study of the embedded topology of plane curves, a curve of the form $C+D$ has been an interesting object as in \cite{shimada2003, shirane2017}.
In this note, we consider an algorithm to produce contact curves to a smooth cubic and apply it to construct Zariski tuples explicitly.
Let $E$ be an elliptic curve given by an affine equation
\[
E: y^2 = f(x),  \, f(x) = x^3 + ax^2 + bx + c.\, \, a, b, c \in \CC.
\]
 Put $O = [0, 1, 0] \in E$ and $\dot{+}$ denotes the addition on $E$ with $O$ as the zero element.
 Let $D$ be a plane curve of degree $r$ such that the divisor $D|_E$ on $E$ defined by $D$ is of the form
 \[
  D|_E = n \left (\sum_{i=1}^d P_i \right ),
 \]
 for some $n > 0$. Put $\fd_D = \sum_{i=1}^dP_i$. Then we have
 $r = nd/3$ by B\'ezout theorem and 
 \[
 n(\fd_D -  dO) \sim 0.
 \]
 Hence $\fd_D - r O$ in $\pic^0(E)$ gives a class whose order divides $n$.

 \begin{defin0}\label{def:n-contact-curve}{\rm 
  We call $D$ a weak $n$-contact
 curve to $E$   if $D|_E = n \left (\sum_{i=1}^d P_i \right ) +  sO$ for some non-negative integer $s$.
 In particular, we call a plane curve $D$  an $n$-contact
 curve to $E$ if $s = 0$. 
 }
 \end{defin0}

 As it can be  seen in \cite{shimada2003, shirane2017},  if $D$ is smooth, plane curves of 
 the form $E+D$ have an  interesting property from a viewpoint of the embedded topology of plane curves.
 Although the existence of such curves was proved in \cite{shimada2003} intrinsically, we here consider a more
 explicit and systematic construction of such curves in this article. In the study of the embedded topology of plane curves, it is one of important steps to construct curves with the same combinatorics having different
 geometric nature and our construction can be another new approach.
 Our main tools are  representations of 
 divisors of degree $0$ on a hyperellitpc curve considered in \cite{leitenberger, mumfordII} using 
 polynomials, the addition on ${\pic}^0(E)$ based on such representations and divisions on
 $\CC[x, y]$.  The representation for a divisor of degree $0$ given in \cite{mumfordII} is
 called the {\it Mumford representation} and in cryptography it plays an important role to describe
 the addition on the Jacobian of a hyperelliptic curve. In this note,
 we first give an method to construct (weak) $n$-contact curves to a cubic based on representations
 given in \cite{leitenberger}
 (Theorem~\ref{thm:n-contact-curve}, Algorithm 1) and
 we apply it to  study the embedded topology of plane algebraic curves.
 This can be considered as  another
 application of such representations for divisors of degree $0$. 
 Our explicit examples for Zariski $N$-ple based on our costruction  are as follows:

 \begin{itemize}
 
  \item A Zariski triple for  cubic and $4$-contact quartic arrangements.
  
   \item A Zariski quartet for  cubic and $6$-contact sextic or cubic and $8$-contact octic arrangements.
   
   \end{itemize}
   
   Our organization of this paper is as follows:
   In \S 1, we give a brief summary on representations for divisors of degree $0$ on a hyperelliptic curve and
   divisions on $\CC[x, y]$. In \S\S 2 and 3, we explain how we construct (weak) $n$-contact curves.
   We give some examples in \S 4  based on our method. In \S 5, we give examples of Zariski tuples as above.
     

\section{Divisors on $E$ and their representation}

\subsection{Semi-reduced divisors for  a hyperelliptic curve}\label{subsec:1-1}

For a smooth projective curve $C$, we denote the group of divisors on $C$
by $\Div(C)$.  The subgroup of $\Div(C)$ generated by divisors of degree $0$ is
denoted by $\Div^0(C)$.

Let $C$ be a hyperelliptic curve of genus $g$ given by 
\[
C: y^2 = x^{2g+1} + c_1x^{2g} + \ldots + c_{2g+1}.
\]
We denote the point at infinity by $O$.
We here introduce  {\it semi-reduced} and {\it reduced divisors} on a hyperelliptic curve used in the study of hyperelliptic cryptosystem
(\cite{cantor, MWZ}). 

\begin{rem}{\rm 
Note that the terminology {\it reduced divisor} in the below is different from usual sense. For example,
as in \cite[Definition, p. 257]{iitaka}, 
$\fd$ is said to be a reduced divisor if
$\fd$ is of the form $\sum_i^d P_i$, $P_i$'s are all distinct. Namely all the non-zero coefficients
in $\fd$ is $1$.  When we use the terminology in this sense, we say {\it a reduced divisor in
usual sense}.
}
\end{rem}

Let $\iota : C \to C$ be the hyperelliptic involution $(x, y) \mapsto (x, -y)$ and let $\pi: C \to \PP^1$ be the
quotient map induced by $\iota$. Put $O_{\PP^1}:= \pi(O)$. $\Div^0(\PP^1)$ is 
generated by divisors of the form $Q - O_{\PP^1}$ $Q \in \PP^1$.  
$\pi^*\Div^0(\PP^1)$  denotes the subgroup of  $\Div^0(C)$ given by
the pull-back  of $\Div^0(\PP^1)$. If we put $\pi^*Q = Q^+ + Q^-$, 
$\pi^*\Div^0(\PP^1)$ is generated by elements of the form $Q^+ + Q^- - 2O$.
For an effective  divisor $\fd$ on $\mcC$ with
$\fd = \sum_{P\in C} m_PP$
by considering points of the form $P + \iota(P)$ contained in 
$\fd$, we have a decomposition 
$\fd = \fd_{\sr} + \fd_o$
such that 
\begin{enumerate}

\item[(i)] $\fd_o \in \pi^*\Div(\PP^1)$, and

\item[(ii)] if we write $\fd_{\sr} = \sum_P m'_P P$, then
       (a) $m'_P = 1$ if $m'_P > 0$ and $P$ is a ramification point of $\pi$, and 
       (b) $m'_{\iota(P)} = 0$ if $m'_P > 0$ and $P$ is not a ramification point of 
       $\pi$.
\end{enumerate}

Based on \cite{galbraith}, we define semi-reduced and reduced divisors as follows:

\begin{defin}\label{def:semireduced}{\rm (\cite{cantor, MWZ})
Let $\fd$ be an effective divisor on $\mcC$. 
\begin{enumerate}

 \item[(i)] $\fd$ is called semi-reduced, if $\fd_o = 0$ and $O$ does not appear in $\fd$.

\item[(ii)] A semi-reduced divisor $\sum_i m_i P_i$ is said to be reduced if $\sum_im_i \le g$.

\end{enumerate}

}
\end{defin}

\begin{rem}{\rm 
 In \cite{cantor, MWZ}, a semi-reduced divisor is defined to be a divisor of degree $0$, i.e.,
a  divisor of the form $\fd - (\deg\fd) O$, where $\fd$ is semi-reduced as above.
}
\end{rem}

For an effective
 divisor $\fd_i = \sum_j m_{i, j}P_j \, \, (i = 1, 2)$, we define $\gcd(\fd_1, \fd_2)$ as follows:
\[
\gcd(\fd_1, \fd_2):= \sum_j\min(m_{1, j}, m_{2, j})P_j.
\]

\subsection{Division algorithms on $\CC[x,y]$ and polynomial functions on hyperelliptic curves}

 Let $C$ be a hyperelliptic curve as in \ref{subsec:1-1} and let 
  $\langle y^2 - f\rangle \subset \CC[x,y]$ be the ideal generated by $y^2-f$.
The quotient ring $\CC[x, y]/\langle y^2 -f \rangle$ is said to be the coordinate
ring of $C$ and we denote it by $\CC[C]$.  The quotient field of $\CC[C]$ is the rational function field $\CC(C)$ of $C$.
 An element of $\CC[C]$ is called
a polynomial function. For $g \in \CC[x, y]$, its class in $\CC[C]$ gives a
polynomial function on $C$, which we denote by $[g]$. 
For two polynomial  functions $[h_1], [h_2]$, we define a divisor $\divi([h_1], [h_2])$ by
\[
\divi([h_1], [h_2]) = \gcd(([h_1])_0, ([h_2])_0) - \deg( \gcd(([h_1])_0, ([h_2])_0) )O,
\]
where $([h_i])_0$ ($i = 1, 2$) denote the zero divisors of $[h_i]$ ($i = 1,2$), respectively.

Define a $\CC[x]$ (resp. $\CC[y]$)-submodule $\amari(y)$ (resp. $\amari(x^{2g+1})$) of $\CC[x, y]$ as follows:

\[
 \amari(y^2)  =  \{b_0 + b_1y \mid b_i \in \CC[x]\}, \quad 
 \amari(x^{2g+1})  =  \left. \left\{\sum_{i=0}^{2g} c_{i}x^{i} \right | c_i \in \CC[y], 0 \le  i \le 2g \right\}.
\]
Then we infer that
$\CC[C] \cong \amari(y^2)$ (resp. $\amari(x^{2g+1})$) as $\CC[x]$ (resp. $\CC[y]$) modules, respectively.
We describe these isomorphisms  based on division algorithm on $\CC[x, y]$.
Let $>_1$ and $>_2$ denote lexicographic orders given by $y>x$ and $x>y$,
respectively. Since $\langle y^2 - f\rangle$ is principal, $y^2- f$ is a Groebner
basis for both $>_1$ and $>_2$.  For $g \in \CC[x, y]$, we consider
divisions with respect to $>_1$ and $>_2$ and two presentations of $g$ are given as  follows:
\begin{eqnarray*}
g & = & q_1(x, y)(y^2  - f) + r_1 \\
& = & q_2(x, y)(y^2 -f) + r_2,
\end{eqnarray*}
where $r_i$ ($i = 1, 2$) denote the remainders with respect to $>_i$ ($i = 1, 2$).
 By \cite[Proposition 1, Ch. II, \S 6]{CLO}, 
 (i) no term of $r_1$ is   divisible
  by  $y^2$, i.e., $r_1 \in \amari(y^2)$ and
  %
  (ii) no term of $r_2$ is   divisible by $x^{2g+1}$, i.e, $r_2 \in \amari(x^{2g})$. 
  Moreover, 
  \begin{enumerate}
  
  \item[(i)] if $g$ is represented as 
  $g =  g_1 + r'_1$ such that $g_1 \in \langle y^2 - f \rangle$  and no term of $r'_1$ is   divisible
  by  $y^2$, then  $r_1 =r'_1$, and
  
  \item[(ii)]  if $g$ is represented as 
  $g =  g_1 + r'_2$ such that  $g_1 \in \langle y^2  - f \rangle$ and no term of $r'_2$ is   divisible
  by  $x^{2g+1}$, then  $r_2 =r'_2$.
  
  \end{enumerate}
  
  We denote ${\mathrm{NF}}_i(g) = r_i$ ($i = 1, 2$). From the above observation,  for $g_1, g_2 \in
  \CC[x, y]$,  we see that ${\mathrm{NF}}_i$ $(i = 1, 2)$ satisfy the following properties:
  
  \begin{enumerate}
   \item[(i)] ${\mathrm{NF}}_i (g_1 + g_2) = {\mathrm{NF}}_i (g_1) + {\mathrm{NF}}_i (g_2)$ $(i = 1, 2)$.

  \item[(ii)] ${\mathrm{NF}}_1 (bg) =b {\mathrm{NF}}_i (g)$ for  $b \in \CC[x], g \in \CC[x, y]$.
  
  \item[(iii)] ${\mathrm{NF}}_2 (cg) = c {\mathrm{NF}}_i (g)$ for  $c \in \CC[y], g \in \CC[x, y]$.
  
  \end{enumerate}
  Hence ${\mathrm {NF}}_1 : \CC[x,y] \to \amari(y^2)$ and ${\mathrm {NF}}_2: \CC[x, y] \to
  \amari(x^{2g+1})$ are $\CC[x]$ and $\CC[y]$ module homomorphisms, respectively.
  Since $[g] = [{\mathrm{NF}}_i(g)]$ and ${\mathrm {NF}}_1(g) = g$ (resp. ${\mathrm {NF}}_2(g) = g$)
  if $g \in \amari(y^2)$ (resp. $g \in \amari(x^{2g+1})$),  the canonical homomorphism $[\bullet] : 
  \CC[x, y] \to \CC[C]$ induces module isomorphisms between $\amari(y^2)$ (resp. $\amari(x^{2g+1}$)) and
  $\CC[C]$, respectively.

\subsection{The case of $g=1$}

In this section, we consider the case when $g=1$.
We keep our notation in the Introduction. Let $\fd$ be an effective divisor.
We first summarize some facts on a method to present $\fd - d O, d = \deg \fd$ by two polynomials
in $\CC[x, y]$ based on \cite{leitenberger}. Let $\iota : E \to E$ be an involution on $E$ given by
$(x, y) \mapsto (x, -y)$, which gives the inversion with respect to  $\dot{+}$.

For $n \in \ZZ_{\ge 0}$, put
\[
\bbL(nO):= \{\varphi \in \CC(E)\setminus \{0\} \mid (\varphi) + nO \ge 0\} \cup \{0\} \subset \CC[E],
\]
where $(\varphi)$ is a principal divisor given by $\varphi$.

Here are some remarks:

\begin{itemize} 

\item $\bbL(nO)$ is a vector space over $\CC$ and $\dim\bbL(nO) = n$. 

For $n \ge 3$,  we have 
\[
\mcB_1: [1], [x], \ldots, [x^{\frac{n - \epsilon}2}], 
 [y],  [xy], \ldots, [x^{\frac{n-4+ \epsilon}2}y], 
\]
where $\epsilon$ is the parity of $n$ as a basis of $\bbL(nO)$. For $\bbL(3nO)  (n \ge 1)$, we give another basis for later use as follows:
\[
\mcB_2: [1], [y], \ldots, [y^n], 
[x], [xy], \ldots, [xy^{n-1}], 
[x^2], [x^2y], \ldots, [x^2y^{n-2}]
\]
For $g \in \CC[x, y]$ such that $[g] \in \bbL(nO)$, 
Note that $\mcB_1$ is a basis corresponding to ${\mathrm {NF}}_1(g)$,
 while $\mcB_2$ is a basis corresponding to  ${\mathrm {NF}}_2(g)$.
 
\item Assume that $\fd \neq \iota^*\fd$.  Since $\dim \bbL((d+1)O) = d+1$ there exists $b \in \CC[E],
\deg_y b= 1$ such that $(g(x, y)) = \fd + P_o - (d+1)O$. Note that $b$ is unique up to non-zero constant.

\item If $\fd = \sum_{i=1}^d P_i$ and $b \in \CC[E]$ with $(b) = \fd + P_o - (d+1)O$, then by the definition of $\dot {+}$ we have
\[
P_1 \dot{+} \dots \dot{+} P_d \dot{+} P_o = O.
\]

 \item For $f_1, f_2 \in \CC[x, y]$ with ${\mathrm{NF}}_1(f_i) = f_i$, $f_1 = f_2$ in $\CC[x, y]$ if
and only if $[f_1] = [f_2]$, since $\CC[E] \cong \amari(y^2)$ as a
$\CC[x]$-module. 

Also for $f_1, f_2 \in \CC[x, y]$ with ${\mathrm{NF}}_2(f_i) = f_i$, $f_1 = f_2$ in $\CC[x, y]$ if
and only if $[f_1] = [f_2]$, since $\CC[E] \cong \amari(x^{2g+1})$ 
as a
$\CC[y]$-module. 

\item For $f_1, f_2 \in \CC[x, y]$ with ${\mathrm{NF}}_i(f_j) = f_j$, $(i, j = 1, 2)$
if $([f_1]) = ([f_2])$ as divisors on 
$E$, then $f_2 = cf_1, \, c \in \CC^{\times}$.

\end{itemize}

By the above remarks and  the representation of $\fd$  based on \cite{leitenberger},  we have the following proposition:

\begin{prop}\label{prop:curve-1}{ For an effective divisor $\fd = \sum_{j=1}^d P_j$
($P_j$'s are not necessarily distinct), there exists  $b_{\fd}(x, y) \in \CC[x, y]$ as follows:

\begin{enumerate}
 \item[\rm{(i)}] $b_{\fd}(x , y) = b_0(x) + b_1(x)y$  for some $b_i(x) \in \CC[x] (i = 0, 1)$, unique up to non-zero
 constant  with
 $\deg b_0 \le \frac{d+1 - \epsilon}2, \deg b_1 \le \frac{d-3 + \epsilon}2$, where $\epsilon$ denotes the parity of $d+1$.
 
 \item[\rm{(ii)}] The  rational function on $E$ given by $b_{\fd}$, which we also denote by $b_{\fd}$, gives a divisor $(b_{\fd}) = \fd + P_o - (d+1)O$. Note that $P_1\dot{+}\ldots \dot{+} P_d 
 \dot{+} P_o =O$.
 
 \item[\rm{(iii)}] Let $D_{b_{\fd}}$ be a plane curve given by $b_{\fd} = 0$. Then the divisor on $E$
 cut out by $D_{b_{\fd}}$ is of the form $\fd + P_o + \nu O$, where $\nu = 0, 1$ or $2$
 satisfying $d + 1 + \nu \equiv 0 \bmod 3$.
 
 \item[\rm{(iv)}] If we denote $P_i = (x_i, y_i)$ $(i =1, \ldots, d)$,  we have
\[
\fd - dO= \divi([\Pi_{i=1}^d(x - x_i)], [b_{\fd}]).
\]
 \end{enumerate}
}
\end{prop}


\section{Construction of weak $n$-contact curves to a smooth cubic}\label{sec:weak-n-contact}

  In what follows, we use $g$ for $[g]$, $g \in \CC[x, y]$ for simplicity.
For other notation, we keep that in \S  1.

 Let $\fd = \sum_{i=1}^d P_i$ such that
$P_{\fd} := P_1\dot{+} \ldots \dot{+} P_d = (x_{\fd}, y_{\fd})$ gives a torsion of order $n$. 
Let $b_{\fd} := b_0(x) + b_1(x) y \in \amari(y^2) \cong \CC[E]$ be an element in $\bbL((d+1)O)$ such that
$(b_{\fd}) = \fd + \iota^*P_{\fd} - (d+1)O$.  Let $\xi := \xi_0(x) + \xi_1(x)y$ be an element in $\CC[E]$ such
that $(\xi) = nP_{\fd} - nO$. Note that $\iota^*\xi  =  \xi_0(x) - \xi_1(x)y = n\iota^*P_{\fd} - nO$. Hence
$\xi\iota^*\xi = c (x - x_{\fd})^n$ for some $c \in \CC^{\times}$. 
Under these circumstances, we have a weak $n$-contact curve  as follows:

\begin{algorithm}
\caption{Weak $n$-contact curve}
\label{alg:weak-n-contact}
\begin{algorithmic}[1]
\REQUIRE the polynomial function $b_{\fd}$
\ENSURE A defining equation $b_{n\fd}$ of a weak $n$-contact curve $D_{nb_{\fd}}$ such that 
$(b_{n\fd}) = n\fd - nd O$ and 
$D_{b_{n\fd}}|_E = n\fd + (d_1- nd)O$.
\STATE Compute ${\mathrm{NF}}_1(b_{\fd}^n\xi)$ and denote it by  $c_0(x) + c_1(x)y$. 
\STATE Put $b_{n\fd} = (c_0(x) + c_1(x)y)/(x - x_o)^n$
\RETURN $b_{n\fd}$.
\end{algorithmic}
\end{algorithm}

\begin{thm}\label{thm:n-contact-curve}{ The  polynomial $b_{n\fd}$  in Algorithm 1 gives a defining equation of
a weak $n$-contact curve such that
$D_{nb_{\fd}}|_E = n\fd + (3d_1 - nd)O$, where $d_1 = \deg b_{n\fd}$.}
\end{thm}

\proof  We first note that for $g \in \CC[x, y]$, $g$ and ${\mathrm {NF}}_1(g)$ defines the same rational
function on $E$. Hence we have the equalities of divisors:
\begin{eqnarray*}
(c_0(x) + c_1(x)y )
& = & ({\mathrm {NF}}_1(b_{\fd}^n\xi)) \\
&= & (b_{\fd}^n\xi) \\
& = & n(\fd - dO) + n(P_{\fd} + \iota^*P_{\fd} - 2O).
\end{eqnarray*}

  Since $P_{\fd} + \iota^*P_{\fd} - 2O \sim 0$, i.e., $n(\fd - dO) \sim 0$, there exists 
  $b_{n\fd} = b_0(x) + b_1(x) y 
  \in \amari(y^2)$, unique up to constant, with $(b_{n\fd}) = n(\fd - dO)$. As
  \[
  (b_{n\fd}(x- x_{\fd})^n) = (c_0(x) + c_1(x)y)
    \]
  and 
  ${\mathrm {NF}}_1(b_{n\fd}(x- x_{\fd})^n) = {\mathrm {NF}}_1(c_0(x) + c_1(x)y)$,
  $b_{n\fd}(x- x_{\fd})^n = c (c_0(x) + c_1(x)y)$ in $\CC[x, y]$  for some $c \in \CC^{\times}$.
  Hence $D_{b_{n\fd}}$ is a weak $n$-contact curve with the desired property as above since 
  $D_{b_{n\fd}}$ may pass through $O$.
\endproof


\section{Construction of $n$-contact curves for semi-reduced divisor of degree divisible by $3$}

We keep our notation as before. Let $D$ be an $n$-contact curve  with
$D|_E = n(\sum_{i=1}^dP_i)$.  Then we have $nd = 3\deg D$.  Hence note that
$nd$ is divisible by $3$.

We put $ \fd = \sum_{i=1}^{d} P_i$, where $P_i$' s are not necessarily
distinct. Let $P_{\fd}$ be the point given by $P_1\dot{+} \ldots  \dot{+} P_{d}$ and we assume that $P_{\fd}$ is 
an $n$-torsion. Put $P_i = (x_i, y_i)$. By Proposition~\ref{prop:curve-1}, we have $b_{\fd} \in \amari(y^2)$ such that
$b_{\fd}(x, y) = b_0(x) + b_1(x) y, \quad \deg b_0 \le \frac {d + 1 - \epsilon}2, \, \, \deg b_1 \le \frac {d - 3 + \epsilon}2$, 
$\epsilon =$ the parity of $d + 1$ and 
$\fd - dO= \divi \left ([\Pi_{i=1}^{d} (x - x_i)], [b_{\fd}]\right )$.
We apply Algorithm 1 to $b_{\fd}^n$ and we denote the resulting element in $\amari(y^2)$ by
$b_{n\fd}$. As a rational function on $E$, the divisor of $b_{n\fd}$ is
\[
(b_{n\fd}) = n\left (\sum_{i=1}^{d} P_i  - d O\right ), \quad  D_{b_{n\fd}}|_E = n\fd + (d_1- nd)O
\]
Let $h_{n\fd} := \mathrm{NF}_2(b_{n\fd}) \in \amari(x^3)$. Since $b_{n\fd}$ and $h_{n\fd}$ define
the same rational function on $E$, we have
\[
(h_{n\fd}) = n \left (\sum_{i=1}^{d} P_i  - d O \right ).
\]
For the basis $\mcB_2$, $y^{nd/3}$ is a unique element which has a pole of order $nd$ at $O$. Hence $h_{n\fd}$ is 
presented by the basis $\mcB_2$ such that the coefficient of $y^{nd/3}$ is not zero. This means that the
curve $D_{h_{n\fd}}$ of degree $nd/3$ given by $h_{n\fd} = 0$ does not pass through $O$, i.e., $D_{h_{n\fd}}|_E =  n(\sum_{i=1}^{d} P_i)$.
Hence $D_{h_{n\fd}}$ is an $n$-contact curve.

\begin{rem}{\rm Note that the plane curve $D_{b_{n\fd}}$ given by $b_{n\fd} = 0$ is a weak $n$-contact
curve since it may intersect $E$ at $O$.
}
\end{rem}

 

\section{Examples}\label{sec:examples}

In this section, we consider some explicit examples for the case of degree $3$. 
Let $T$ be a torsion of order $n$.  Let $\fd$ be a semi-reduced divisor of degree $3$ such
that $P_{\fd} = T$.
In this case,
we can choose $b_{\fd}$ in such a way that  $b_{\fd} = y + b_0$ $\deg b_0 \le 2$. Hence
the pair $(\prod_{i=1}^3(x - x_i), -b_0)$ is the Mumford representation of $\fd$.
Thoughout this section,  $\xi$ denotes an element in $\amari(y^2)$ 
such that $(\xi) = n(T - O)$ as \S \ref{sec:weak-n-contact}.  We make use of Maple 19 for our computation.

\subsection{$2$-torsions}\label{subsec:2torsion}
Let $E_2$ be an elliptic curve with a $2$-torsion.
A $2$-torsion  is given by $(x_o, 0)$. Here $x_o$ satisfies $f(x_o) = 0$, i.e., $E_2$ is given by
$E_2: y^2 - f_2, f_2(x) = (x - x_o)(x^2 + c_1x + c_2)$ and $\xi$ is given
by $x - x_o$.  A rational function $b_{\fd}$ with $(b_{\fd}) = \fd + \iota^*T - 4O$ is of the
form $y + a(x - x_o)(x-b)$. We apply our method  in \S 3 to this case and $h_{2\fd}$ is given by
\begin{eqnarray*}
h_{2\fd} & = &
(-2a^2b-a^2c_1+1)x^2+2axy+(a^2b^2+2a^2bx_o+a^2c_1x_o-a^2c_2+c_1)x \\
& & +a^2y^2-2aby-a^2b^2x_o+
a^2c_2x_o+c_2
\end{eqnarray*}

\subsection{$3$-torsions}\label{subsec:3torsion}
Let $E_3$ be an elliptic curve with a $3$-torsion
Let $T = (x_o, y_o)$ be a $3$-torsion of $E_3$ and let $L_T: y = mx + n$ be the tangent line at $T$.  Then
 $E_3$ is given by 
$y ^2 = f_3, f_3(x) = (x - x_o)^3 + (mx+ n)^2$.
The equation of $L_T$ gives $\xi$ as $(mx  +n -y) = 3(T - O)$.  A  rational 
function $b_{\fd}$ with $(b_{\fd}) = \fd + \iota^*T - 4O$ is of the
form $y - r(x - x_o)(x-s) + mx_o +n$. We apply our method  in \S 3 to this case and we obtain
$b_{3\fd}$ and $h_{3\fd}$.  For simplicity, we consider the case when $m = n = r = 1$.\\
In this case, $(b_{\fd}) = \fd + \iota^*T - 4O$ is of the
form $y - (x - x_o)(x-s) + x_o +1$. We apply Algorithm 1 and have

\begin{eqnarray*}
b_{3\fd} & = & \mathrm{NF}_1(b^3\xi)/(x - x_o)^3 \\
&=& s^3x + s^3y - 6s^2xx_o - 3s^2xy + 3s^2x_o^2 - 3sx^3 + 12sx^2x_o + 3sx^2y \\
&& - 6sxx_o^2 + 2x^4 - 6x^3x_o - x^3y + 3x^2x_o^2 + s^3 + 3s^2y - 9sxx_o- 3sxy + 6sx_0^2 \\
&& - 3sx_oy + 6x^2x_o - 3xx_o^2 + 3xx_oy - x_o^3 + 3s^2 - 6sx + 3sx_o + 3sy + 5x^2 \\
&& - 6xx_o - 2xy + 3x_o^2 - 3x_oy + 3s - x + 3x_o - y + 3
\end{eqnarray*}
we have $h_{3\fd}=\mathrm{NF}_2(b_{3\fd})$ as follows:
\begin{eqnarray*}
h_{3\fd}& = &x^2y(3s - 3x_o + 1) + x^2(3sx_o - 3x_o^2 + 3s + 3) + 2xy^2  + xy(-3s^2 + 3x_o^2 \\
&&  - 3s + 3x_o) + x(s^3 - 6s^2x_o + 3sx_o^2 + 2x_o^3 - 9sx_o + 3x_o^2 - 6x_o + 1)  - y^3 \\
&& + y^2(-3s - 2)  + y(s^3 - x_o^3 + 3s^2 - 3sx_o + 3s - 3x_o) + 3s^2x_o^2 - 3sx_o^3 + s^3  \\
&& + 6sx_o^2  - 3x_o^3 + 3s^2 + 3x_os + 3x_o^2 + 6s + 3x_o + 5
\end{eqnarray*}
We can check that a curve $C_3$ given by $h_{3\fd}=0$ is smooth for general  $s$.

\subsection{$4$-torsions}\label{subsec:4torsion}
Let $E_4$ be an elliptic curve with a $4$-torsion.
 Let $T=(x_o, y_o)$ be a $4$-torsion and let $P:= [2]T$. Note that $\iota^*T = (x_o, - y_o)$. As $P$ is a $2$-torsion, we choose 
affine coordinates so that $P = \iota^*(P) =  (0,0)$. Then we may assume that $E_4$ is given by
$E_4: y^2 - f_4,  f_4(x) = x(x^2 + c_1x + c_2)$ 
and the tangent line $L_T$ at $T$ is given by
$\displaystyle{y = mx, \quad m =  f'(x_o)/2y_o}$.
As $(y - mx) = 2T + \iota^*P - 3O$, we have
$((y-mx)^2) = 4(T -O) + 2(T - O).$
Hence $\mathrm{NF}_1((y -mx)^2) = x(x^2 + (c_1 + m)x + c_2 - 2my)$ and 
 $\xi$ is given by $\xi = x^2 + (c_1 +m)x + c_2 - 2my$.
Note that we make use of Algorithm 1 to obtain $\xi$.
A  rational 
function $b_{\fd}$ with $(b_{\fd}) = \fd + \iota^*T - 4O$ is of the
form $y - r(x - x_o)(x-s) + mx_o$. 
We apply our method  in \S 3 to this case and we obtain
$b_{4\fd}$ and $h_{4\fd}$.   
 
  By putting  $c_1 = -(2t-1), c_2 = t^2,
t \in \CC\setminus \{0, 1/4\}$,  we have an elliptic curve $E_4$ with a parameter $t$.
This example given in \cite[p. 57]{silverman-tate} . In this case, $T  = (t, t), P = (0, 0)$ and $m = 1$.
For $b_{\fd}$,  we consider the case of $r = 1$ and $s =2$ for simplicity.

\begin{eqnarray*}
b_{4\fd} & = & \mathrm{NF}_1(b^4 \xi)/(x - t)^4 \\
& = &(t^2x^4 - 2tx^5 + x^6 + 2t^2x^3 - 8tx^4 + 4tx^3y + 8x^5 - 6x^4y - 5t^2x^2 + 26tx^3\\
&&  - 8tx^2y - 29x^4 + 16x^3y + 2t^2x - 20tx^2 + 8txy + 32x^3 - 18x^2y + t^2 + 2tx - 4ty\\
&&  - 11x^2 + 12xy + 2x - 2y)
\end{eqnarray*}
We have $h_{4\fd}=\mathrm{NF}_2(b_{4\fd})$ as follows:
\begin{eqnarray*}
 h_{4\fd}&=&x^2y^2(2t + 6) + x^2y(-10t^2 + 44t - 40) + x^2(8t^3 - 80t^2 + 168t - 79) - 6xy^3+ xy^2(-t^2\\
 && + 6t - 36) + xy(8t^3 - 22t^2 + 8t + 12) + x(-7t^4 + 52t^3 - 66t^2 + 2t + 2) + y^4 + y^3(-8t\\
 && + 22) + y^2(7t^2 - 52t + 68) + y(-4t - 2) + t^2
\end{eqnarray*}
In the case, the curve given by $h_{4\fd}=0$ has a singular point at $[1, 0, 0]$($x = X/Z, y = Y/Z$).  In oder to obtain a smooth $4$ contact curve,
 we consider a curve $C_4$ given $\tilde{h}_{4\fd}:=h_{4\fd}+(x+y+1)(y^2-f_4)= 0$.

\subsection{$6$-torsions}

By \cite[p.238, 8.13]{silverman}, we have an
 elliptic curve $E_6$ ($ t \neq 1, 2, 10/9$) with 6-torsion as follows:
\[
E_6 :  y^2  = f_6(t, x),  \,\,  f_6(t, x)  =   x^3-\left (\frac 34 t^2-3t+2 \right )x^2+\frac12(-t^2+3t-2)tx+\frac14(-t^2+3t-2)^2 
\]

$E_6$ has torsion points for general $t$as follows:
\begin{center}
\begin{tabular}{|c|c|} \hline
 $6$-torsion & 
$\pm T = [0, \pm \frac12 (t-1)(t-2)]$  \rule[-1.5mm]{0mm}{5mm}\\ \hline
$3$-torsion &
$\pm [2]T =\left [t^2-3t+2,  \pm\frac 12 (t-2)(t-1)^2 \right ]$ \rule[-1.5mm]{0mm}{5mm}\\ \hline
$2$-torsion & $[3]T =[-3]T =[-t + 1, 0]$ \\ \hline
\end{tabular}
\end{center}

Here are some tips to compute an equation of $6$-contact curve $D$ to $E$. 
Basically our  $4$-torsion-case method works, and it is a key to find $\xi$, i.e., a rational function which gives
a divisor $6T - 6O$.  More precisely it is as follows:

\begin{enumerate}

\item[Step. i] Find the equation for the tangent line $L$ at $T$, and denote it by $l$. Note that $(L|_E) = 2T + \iota^*[2]T - 3O$.

\item[Step. ii] $l^3$ gives a rational function so that $(l^3) = 6T + 3\iota^*[2]T - 9O = 6(T - O) + 3(\iota^*[2]T - O)$. This
means that $\NF_1(l^3)$ also gives a rational function with same divisor. Let $L_{[2]T}$ be the tangent line
at $[2]T$ and let $l_1$ be its defining equation, which gives a divisor $3([2]T - O)$.
Therefore $\NF_1(l^3l_1)$ is gives a divisor $6(T - O) + 3(\iota^*[2]T + [2]T - 2O)$.
\item[Step. iii] Now since $\iota^*[2]T + [2]T - 2O$ is a divisor of $x - (t^2-3t+2)$ ,  $\xi = \NF_1(l^3l_1)/(x - (t^2-3t+2))^3$:
\begin{eqnarray*}
\xi = t^3 - 3t^2x + 2x^3 - 5t^2 + 8tx - 2ty + 4x^2 + 4xy + 8t - 4x + 4y - 4
\end{eqnarray*}
%

\end{enumerate}
Let $T = (x_o, y_o)=[0, \frac12 (t-1)(t-2)]$ be a $6$-torsion of $E$,
a rational function $b_{\fd}$ with $(b_{\fd}) = \fd + \iota^*T - 6O$ is of the
form $y - r(x - x_o)(x-s) + y_o$.  For simplicity, We consider the case when $r=1$, $s=4$ and $t=3$. 

\begin{eqnarray*}
b_{6\fd} & = & 64\mathrm{NF}_1(b_{\fd}^6 \xi)/x^6\\
& = &128x^9-2432x^8-512x^7y+24864x^7+8832x^6y-173184x^6-57024x^5y+738248x^5\\
&& +107168x^4y-1310712x^4+590592x^3y-1918138x^3-3714312x^2y+11061932x^2\\
&& + 7111844xy-12378399x-3541074y+3545170.
\end{eqnarray*}
We have $h_{6\fd} = \mathrm{NF}_2(b_{6\fd})$ as follows:
\begin{eqnarray*}
h_{6\fd}& = & \frac{68872271}{32}-2528x^2y^4-512xy^5+128y^6-64608x^2y^3+27256xy^4+9088y^5 \\
&& +\frac{1997169}{2}y^2x^2 +177536y^3x-201632y^4-3294680x^2y-\frac{21965953}{8}y^2x+352704y^3 \\
&& +\frac{448707487}{128}x^2+8047460xy+\frac{51021297}{32}y^2-\frac{194729737}{32}x-3902866y.
\end{eqnarray*}
In the case, $h_{6\fd}=0$ has a singular point at $[1, 0, 0]$ ($x = X/Z, y = Y/Z$).  We then consider 
$\tilde{h}_{6\fd}:=h_{6\fd}+(x^3+y^3+1)(y^2-f_6(3, x))$.
We can check that a curve $C_6$ given by $\tilde{h}_{6\fd}(x)=0$ is smooth.

\subsection{$8$-torsions}

Again by \cite[p.238, 8.13]{silverman}, we have we have an
 elliptic curve $E_8$ ($t \neq 0, 1/1/2, (2\pm \sqrt{2})/2$) with 8-torsion as follows:
\[
E_8 :  y^2 = f_8(t, x), \,\, f_8(t,x) =  (x - t^4 +t^3)\left (x^2-(2t^3-4t^2+2t-1/4)x-t^6+2t^5-\frac 54 t^4 + \frac 14t^3 \right ) 
\]
For general $t$, $E$ has torsion points as follows:

\begin{center}
\begin{tabular}{|c|c|} \hline
 $8$-torsion & 
$\pm T  = \left [0, \pm(-t^5 + \frac 32 t^4 - \frac 12 t^3) \right ]$   \rule[-1.5mm]{0mm}{5mm}\\ \hline
$4$-torsion &
$\pm [2]T = \left [ t^2(2t-1)(t-1), 2(t-1)^2t^2\left (t- \frac 12\right )^2 \right ]$ \rule[-2.5mm]{0mm}{7mm}\\ \hline
$2$-torsion & $[t^3(t-1), 0]$ \\ \hline
\end{tabular}
\end{center}

%

In order to obtain $\xi$, we repeat a similar argument to that in the case of $6$-torsion, which is as follows:
\begin{enumerate}

\item[Step. i] Find the equation for the tangent line $L$ at $T$, and denote it by $l$. Note that $(L|_E) = 2T + \iota^*[2]T - 3O$.

\item[Step. ii] $l^4$ gives a rational function so that $(l^4) = 8T + 4\iota^*[2]T - 12O = 8(T - O) + 4(\iota^*[2]T - O)$. This
means that $\NF_1(l^4)$ also gives a rational function with same divisor. Let $L_{[2]T}$ be the tangent line
at $[2]T$ and let $l_1$ be its defining equation. Then $2[2]T + [4]T - 3O $.
Therefore $\NF_1(l^4l_1^2)$ is gives a divisor $8(T - O) + 4([2]T  +\iota^*[2]T - 2O) + 2([4]T - O)$.
\item[Step. iii] Now since $\iota^*[2]T + [2]T - 2O$ is a divisor of $x -t^2(2t-1)$,  $\xi = \NF_1(l^4l_1^2)/(x - t^2(2t-1)(t-1))^4(x-t^3(t-1))$
\begin{eqnarray*}
\xi &:=&  4t^{12}-8t^{11}+5t^{10}+16t^9x-t^9-32t^8x+18t^7x+4t^7y+16t^6x^2-3t^6x-2t^6y\\
&& -40t^5x^2+21t^4x^2+12t^4xy-3t^3x^2-4t^3xy-16t^2x^3+8tx^3+8tx^2y-2x^4-x^3\\
&& -2x^2y
\end{eqnarray*}


\end{enumerate}

Let $T = (x_o, y_o)=\left [0, -t^5 + \frac 32 t^4 - \frac 12 t^3 \right ]$ be a $8$-torsion of $E_8$,
a rational function $b_{\fd}$ with $(b_{\fd}) = \fd + \iota^*T - 8O$ is of the
form $y - r(x - x_o)(x-s) + y_o$. 
We consider the case when $r = 1$,  $s=1$ and $t = -1$ for simplicity. Then we have
\begin{eqnarray*}
b_{8\fd} & = & 256\mathrm{NF}_1(b_{\fd}^8 \xi)/x^8 \\
& = &512x^{12}-3840x^{11}-1536x^{10}y+53760x^{10}-4096x^9y-78848x^9\\
&& +226304x^8y-4409664x^8-2783232x^7y+48818400x^7+20539072x^6y\\
&& -283514336x^6-97078784x^5y+1066337424x^5+304250784x^4y\\
&& -2692293822x^4-639609984x^3y+4513159593x^3+874149354x^2y\\
&& -4813270128x^2-702018576xy+2958279813x+250317702y-798728850
\end{eqnarray*}
We have $h_{8\fd} = \mathrm{NF}_2(b_{8\fd})$ as follows:
\begin{eqnarray*}
h_{8\fd}& = & \frac{-986999682916161}{512}-16640x^2y^6-1536xy^7+512y^8-153664x^2y^5+288768xy^6\\
&& +24704y^7 +54601398x^2y^4-2994208xy^5-4303600y^6-974116559x^2y^3-272142150xy^4 \\
&& +69041954y^5 -\frac{157122991069}{16}x^2y^2+\frac{45769674189}{8}xy^3+\frac{11197406491}{8}y^4\\
&& +\frac{41898589232589}{128}yx^2 +\frac{755425355935}{128}y^2x-\frac{1214565010245}{32}y^3-\frac{3844005564383585}{2048}x^2 \\
&&-\frac{26329075721469}{32}xy +\frac{103349774295737}{512}y^2+\frac{2415502623658065}{512}x\\
&& +\frac{10759562218989}{32}y
\end{eqnarray*}
In the case, $h_{8\fd}=0$ has a singular point at $[1, 0, 0]$ ($x = X/Z, y = Y/Z$).  We then consider 
$\tilde{h}_{8\fd}:=h_{8\fd}+(x^5+y^5+1)(y^2-f_8(-1,x))$.
We can check that a curve $C_8$ given by $\tilde{h}_{8\fd}=0$ does not have a singular point.



\section{Zariski tuples}

As an application of our construction of a $n$-contact curve to $E$, we consider Zariski pairs and tuples for a smooth cubic and its $n$-contact curves explicitly.
For definition and terminologies for Zariski pairs or tuples, we refer to \cite{survey} and use terminologies there freely. The combinatorics considered here is
\[
B = E + C,
\]
where $E$ is a smooth cubic and $C$ is a smooth $n$-contact curve with $C\cap E = \{P_1, P_2, P_3\}$ (Note that $\deg C = n$). Such examples have
already been
considered in \cite{shimada2003, shirane2017} and the existence of such curves are given in an intrinsic way. In this section, we construct such examples
for $n = 4, 6, 8$, explicitly.
Our key tool to  distinguish the embedded topology of $B$ is the splitting number of $E$ with respect to a cyclic cover $\phi_C : S_C \to \PP^2$ of degree
$n$ branched along $C$ (see \cite[Corollary 1.4]{shirane2017}). By \cite[Proposition 2.5]{shirane2017}, in order to compute the splitting
number $s_{\phi_C}(E)$, we only need to know the order of $P_1 + P_2 + P_3 - 3O$ in  $\pic^0(E)$, i.e., the order of $T = P_1\dot{+}P_2\dot{+}P_3$.
Let $E$ and $Q_d$ be a smooth cubic and a general smooth curve. 
Let $E_d$ and $C_d$ be elliptic curves and smooth $d$-contact curves given in \S \ref{sec:examples}. For $n \ge 4$ with $d|n, d \ge 2$, consider a pencil
\[
\Lambda_d = |\lambda(n/d C_d) + \mu(E_d + Q_{n-3})|_{[\lambda, \mu] \in \PP^1}.
\]
A general member $D \in \Lambda_d$ is a smooth $n$-contact curve to $E_d$.  Let $\phi : S_{D} \to \PP^2$ be an $n$-cyclic cover of degree $n$.
Then,  by \cite[Proposition 2.5]{shirane2017}, $s_{\phi}(E) = n/d$.

\subsection{The case of $n = 4$} Let us consider the elliptic curve considered in \S\ref{subsec:4torsion}

Let $\Lambda_1:=  |\lambda(4 L) + \mu(E + Q_1)|_{[\lambda, \mu] \in \PP^1}$ and 
$\Lambda_2:= \lambda (2C_2) + \mu (E_2 + Q_2)|_{[\lambda, \mu] \in \PP^1}$. Let $D_i$ $(i = 1, 2)$ be general members of $D_i$ $(i =1,2)$, respectively.
Let  $D_3 = C_4$ in \S\ref{sec:examples}
 
 Let $\phi_i : S_{D_i} \to \PP^2$ be $4$-fold cyclic covers branched along $D_i$ ($i = 1, 2,3$), respectively.
 Then splitting numbers $s_{\phi_i}(E)$ are
$ s_{\phi_1}(E) = 4, s_{\phi_2}(E_2) = 2$, and $s_{\phi_3}(E_4) = 1$.
 Hence $(E + D_1, E+ D_2, E+D_3)$ is a Zariski triple.

\subsection{The case of $n = 6$}
Let $\Lambda_1:=  |\lambda(6 L) + \mu(E + Q_3)|_{[\lambda, \mu] \in \PP^1}$,
$\Lambda_2:= \lambda (3C_2) + \mu (E_2 + Q_3)|_{[\lambda, \mu] \in \PP^1}$ and
$\Lambda_3:= \lambda (2C_3) + \mu (E_2 + Q_3)|_{[\lambda, \mu] \in \PP^1}$ 
 Let $D_i$ $(i = 1, 2, 3)$ be general members of $D_i$ $(i =1,2, 3)$, respectively
Let  $D_4 = C_{6\fd}$ in \S\ref{sec:examples}.
 Let $\phi_i : S_{D_i} \to \PP^2$ be $6$-fold cyclic covers branched along $D_i$ ($i = 1, 2,3, 4$), respectively.
 Then splitting numbers $s_{\phi_i}(E)$ are
$ s_{\phi_1}(E) = 6$,  $s_{\phi_2}(E_2) = 3$, $s_{\phi_3}(E_4) = 2$ and   $s_{\phi_4}(E_6) = 1$.
 Hence $(E + D_1, E+ D_2, E+D_3, E+D_4 )$ is a Zariski quartet.

\subsection{The case of $n = 8$}
Let $\Lambda_1:=  |\lambda(8 L) + \mu(E + Q_4)|_{[\lambda, \mu] \in \PP^1}$,
$\Lambda_2:= \lambda (4C_2) + \mu (E_2 + Q_4)|_{[\lambda, \mu] \in \PP^1}$ and
$\Lambda_3:= \lambda (2C_4) + \mu (E_4 + Q_4)|_{[\lambda, \mu] \in \PP^1}$ 
 Let $D_i$ $(i = 1, 2, 3)$ be general members of $D_i$ $(i =1,2, 3)$, respectively
Let  $D_4 = C_{8}$ in \S\ref{sec:examples}.
 Let $\phi_i : S_{D_i} \to \PP^2$ be $8$-fold cyclic covers branched along $D_i$ ($i = 1, 2,3, 4$), respectively.
 Then splitting numbers $s_{\phi_i}(E)$ are
$ s_{\phi_1}(E) = 8$,  $s_{\phi_2}(E_2) = 4$, $s_{\phi_3}(E_4) = 2$ and   $s_{\phi_4}(E_8) = 1$.
 Hence $(E + D_1, E+ D_2, E+D_3, E+D_4 )$ is a Zariski quartet.

\noindent Ai Takahashi and Hiro-o Tokunaga\\
Department of Mathematical  Sciences, Graduate School of Science, \\
Tokyo Metropolitan University, 1-1 Minami-Ohsawa, Hachiohji 192-0397 JAPAN \\
{\tt tokunaga@tmu.ac.jp}

\end{document}